\scrollmode
\documentclass[12pt]{amsart}
\usepackage{verbatim}
\usepackage{eucal}
\usepackage{amssymb}
\usepackage{mathrsfs}
\usepackage{graphicx}
\usepackage{psfrag}
\newif \ifwide
\newif \ifavnermargin
\def \makemargins{
\ifwide
	\oddsidemargin .25in
	\evensidemargin .25in
	\textwidth 6.00in
\else
\fi
\ifavnermargin
	\headheight=7pt
	\textheight=574pt
	\textwidth=432pt
	\topmargin=14pt
	\oddsidemargin=18pt
	\evensidemargin=18pt
\else	
\fi
}

\widetrue
\makemargins

\theoremstyle{plain}
\newtheorem{theorem}[subsection]{Theorem}
\newtheorem{proposition}[subsection]{Proposition}
\newtheorem{lemma}[subsection]{Lemma}
\newtheorem{corollary}[subsection]{Corollary}

\theoremstyle{definition}
\newtheorem{definition}[subsection]{Definition}
\newtheorem{example}[subsection]{Example}
\newtheorem{exam}[subsection]{Example}

\theoremstyle{remark}
\newtheorem{remark}[subsection]{Remark}

\newcount\timehh\newcount\timemm
\timehh=\time 
\divide\timehh by 60 \timemm=\time
\newcommand{\draftauthor}[1]{\author{#1
    {
      --- \protect \protect\sc\today\ ---
      \ifnum\timehh<10 0\fi\number\timehh\,:\,\ifnum\timemm<10 0\fi\number\timemm
      \protect \, \, \protect \bf DRAFT
    }
  }
}
\newcommand{\RR}{{\mathbb R}}
\newcommand{\CC}{{\mathbb C}}
\newcommand{\ZZ}{{\mathbb Z}}
\newcommand{\QQ}{{\mathbb Q}}
\newcommand{\Proj}{{\mathbb P}}
\newcommand{\ee}{{\mathrm{e}}}
\newcommand{\ii}{{\mathrm{i}}}
\newcommand{\RRR}{{\mathscr{R}}}
\newcommand{\CCC}{{\mathscr{C}}}
\newcommand{\OOO}{{\mathscr{O}}}
\newcommand{\HHH}{{\mathfrak{H}}}
\newcommand{\TTT}{{\mathscr{T}}}
\newcommand{\III}{{\mathscr{I}}}
\newcommand{\MMM}{{\mathscr{M}}}
\newcommand{\ac}[1]{{\rm a.c.}\Bigl(#1\Bigr)}
\renewcommand{\theta}{\vartheta}
\renewcommand{\tilde}{\widetilde}

\renewcommand{\mod}{\bmod}

\DeclareMathOperator{\Tor}{Tor}
\DeclareMathOperator{\Hom}{Hom}

\DeclareMathOperator{\Spec}{Spec}

\DeclareMathOperator{\relint}{rel.int.}
\DeclareMathOperator{\Rank}{rank}
\DeclareMathOperator{\Dim}{dim}
\DeclareMathOperator{\codim}{codim}

\DeclareMathOperator{\Trace}{Trace}
\DeclareMathOperator{\Td}{Td}
\DeclareMathOperator{\Sym}{Sym}
\DeclareMathOperator{\chernchar}{ch}

\begin{document}

\title{Toric Varieties and Modular Forms}

\newif \ifdraft
\def \makeauthor{
\ifdraft
	\draftauthor{Lev A. Borisov and Paul E. Gunnells}
\else
\author{Lev A. Borisov}
\address{Department of Mathematics\\
Columbia University\\
New York, NY  10027}
\email{lborisov@math.columbia.edu}

\author{Paul E. Gunnells}
\address{Department of Mathematics\\
Columbia University\\
New York, NY  10027}
\email{gunnells@math.columbia.edu}

\fi
}

\thanks{Both authors were partially supported by a Columbia
University Faculty Research grant}

\draftfalse
\makeauthor

\ifdraft
	\date{\today}
\else
	\date{August 26, 1999}
\fi


%
%

\begin{abstract}
Let $N\subset \RR^{r}$ be a lattice, and let $\deg\colon N \rightarrow
\CC$ be a piecewise-linear function that is linear on the cones of a
complete rational polyhedral fan.  Under certain conditions on $\deg$,
the data $(N,\deg)$ determines a function $f\colon {\HHH}\rightarrow
\CC$ that is a holomorphic modular form of weight $r$ for the
congruence subgroup $\Gamma_{1} (l) $.

Moreover, by considering all possible pairs $(N ,\deg)$, we obtain a
natural subring ${\TTT} (l)$ of modular forms with respect to
$\Gamma_{1} (l) $.  We construct an explicit set of generators for
$\TTT (l)$, and show that ${\TTT} (l)$ is stable under the action of
the Hecke operators.  Finally, we relate ${\TTT} (l)$ to the
Hirzebruch elliptic genera that are modular with respect to
$\Gamma_{1} (l) $.
\end{abstract}
\maketitle

%
%
\section{Introduction}\label{introduction}
\subsection{}
The construction of arithmetically distinguished automorphic forms
from theta series has a long and rich history.  An early spectacular
manifestation is the computation of the number of representations of
an integer by a sum of four squares, in which theta series are
compared with Eisenstein series (cf. \cite{mumford}).  Another is
Hecke's basis problem, which asks if all modular forms can be
expressed in terms of theta functions of quaternary quadratic forms
associated to orders in definite rational quaternion algebras
\cite{eich1, eich2, hps}.  In \cite{wald}, Waldspurger has made a deep
study of the generation of modular forms by theta series.

In this paper, we construct a subspace $\TTT (l)$ of the level $l$
holomorphic modular forms with character using toric geometry and the
combinatorics of rational polyhedral fans.  We show that our modular
forms are related to products of logarithmic derivatives of theta
functions with characteristic.  Our main results say that $\TTT (l)$
is a natural subspace:
\begin{itemize}
\item It is a finitely generated ring over $\CC$.
\item It is stable under the action of the Hecke operators.
\item It is stable under the Fricke involution.
\item It is stable under Atkin-Lehner lifting.
\end{itemize}
For $l\geq 5$, the generators of $\TTT (l)$ are the
\emph{Hecke-Eisenstein} forms (see Remark~\ref{lang.rmk}).
\subsection{}
To describe our construction in more detail, we must fix notation.
Let $l$ be a positive integer, and let $\Gamma _{1} (l)\subset SL_{2}
(\ZZ )$ be congruence subgroup of matrices satisfying
\[
\Bigl (\begin{array}{cc}
a&b\\
c&d
\end{array}\Bigr)=
\Bigl (\begin{array}{cc}
1&*\\
0&1
\end{array}\Bigr)\mod l.
\] 
For any positive integer $r$, let $\MMM _{r} (\Gamma _{1} (l))$ be the
$\CC$-vector space of weight $r$ holomorphic modular forms for $\Gamma
_{1} (l)$. 

Let $X$ be a complete, possibly singular, toric variety of
dimension $r$, and let $\Sigma \subset \RR ^{r}$ be the corresponding
fan, which we assume to be rational with respect to a lattice $N$.
Let $\deg \colon N \rightarrow \CC$ be a piecewise-linear function
that is linear on the cones of $\Sigma $.  Under certain conditions on
$\deg$, the data $(N,\deg)$ determines a function $f_{N,\deg}\colon
\HHH \rightarrow \CC $ that lies in $\MMM _{r} (\Gamma _{1} (l))$.  We
call $f_{N,\deg}$ a \emph{toric form}.

Essentially, $f_{N,\deg}$ is constructed as the alternating sum over
the cones of $\Sigma $ of a collection of $q$-expansions
\begin{equation}\label{informal}
\sum_{m\in M}\sum _{n\in C\cap N} q^{m\cdot n}\ee^{2\pi \ii \deg (n)},
\end{equation}
where $M$ is the dual lattice of $N$, and $m\cdot n$ is the natural
pairing.  However, we emphasize that \eqref{informal} is a formal
expression, and some effort is required to see that \eqref{informal},
and the sum over the cones of $\Sigma $, has meaning.

\subsection{}
In a sense, our $q$-expansions differ dramatically from those
constructed from quadratic forms, in that they are
``piecewise-linear'' theta series: instead of counting the number of
lattice points on an ellipsoid, we essentially count the number of
lattice points on a polytope, with twisting by roots of unity provided
by $\deg$.  Because of this, one might expect that toric forms have
little arithmetic content, or at best that their span contains only
Eisenstein series.  However, this is definitely not the case.  For
weight $2$ and prime level $l<37$, for example, we have $\TTT_{2}(l) =
\MMM _{2} (l)$.  In general we do not know how much of $\MMM (\Gamma
_{1} (l))$ is captured by $\TTT (l)$, nor can we provide an arithmetic
characterization of the toric forms (however, see
Remark~\ref{conjecture}).
 
\subsection{}
Here is an overview of this paper.  In \S\ref{degree.fn.sect}, we
define toric forms, and show that they are well-defined by
relating them to a construction in homological algebra.  In
\S\ref{chern.sect}, we describe toric forms in terms of theta
functions by interpreting the former as the graded Euler
characteristic of a certain infinite-dimensional bundle over $X$.  In
\S\ref{toric.ell}, we show that toric forms are modular, and give a
set of generators for $\TTT (l)$.  In \S\ref{hecke.sect}, we prove our
main results about the compatibility of $\TTT (l)$ with Hecke theory
and lifting.  Finally, in an appendix (\S\ref{appendix}), we present
background on toric varieties for readers not familiar with them.

%
%
\section{Degree functions and functions on the upper
halfplane}\label{degree.fn.sect}
\subsection{}
Let $N$ be a lattice of rank $r$, and let $\deg\colon N \rightarrow
\CC $ be a piecewise-linear function.  Then $\deg $ is a \emph{degree
function} if there exists a complete rational polyhedral fan $\Sigma $
such that $\deg $ is linear on all cones of $\Sigma$.

Let $M=\Hom _{\ZZ }(N,\ZZ)$ be the dual lattice, and let $m\cdot n$ be
the natural pairing.  Let $M_{\CC} = M\otimes \CC $.  For every cone
$C\in\Sigma$, we can define a map
$$
f_C\colon {\HHH}\times M_\CC \rightarrow \CC
$$
as follows.  For $\tau \in \HHH $, let $q = \ee ^{2\pi \ii }$.  If
$m\in M$ satisfies
\begin{equation}\label{cond.on.m}
m\cdot (C\smallsetminus\{0\})>0,
\end{equation}
then we set 
\[
f_{C} (q, m) := \sum_{n\in C} q^{m\cdot n} \ee^{2\pi\ii \deg(n)},
\]
for all $\tau $ with sufficiently large imaginary part.
For all other $m$ and $q$, we define $f_{C}$ by analytic continuation.  It is
easy to see that $f_C$ is a rational function of $q^{m\cdot n_i}$,
where $\{n_i\}$ is any basis of $N$.  To emphasize the origin of
$f_{C}$ as a sum, we will usually write
$$
\ac{\sum_{n\in C} q^{m\cdot n} \ee^{2\pi\ii \deg(n)}}
$$
instead of $f_C(q,m)$.

\begin{definition}
Let $\deg \colon N \rightarrow \CC $ be a degree function with respect
to the fan $\Sigma $.  Assume that $\deg (d) \not \in \ZZ $ for the
generator of any one-dimensional cone of $\Sigma$.  Then the
\emph{toric form} associated to $(N,\deg)$ is the function
$f_{N,\deg}\colon \HHH\rightarrow \CC $ defined by
$$
f_{N,\deg}(q) := \sum_{m\in M}\Bigl(\sum_{C\in\Sigma}(-1)^{\codim C}
\ac{\sum_{n\in C} q^{m\cdot n} \ee^{2\pi\ii \deg(n)}}\Bigr).
$$
\label{modularfromdegree}
We denote by $\TTT$ the $\CC $-vector space generated by the
toric forms, and by $\TTT _{r}\subset \TTT $ the subspace generated by
those $f_{N,\deg}$ with $\Rank N = r$.
\end{definition}

\begin{exam}
Suppose that $\deg(d)=1/2$ for all generators $d$ of one-dimensional
cones of $\Sigma$, and suppose that the toric variety $X$ associated
to $\Sigma $ is nonsingular.  Then the function $f_{N,\deg}$ is a
normalized elliptic genus of $X$ \cite{borlibg}.  This example was our
major motivation and the starting point for this paper.
\end{exam}

\begin{example}
Let $N = \ZZ ^{2}\subset \RR ^{2}$, and let $\Sigma $ be the fan in
Figure~\ref{p2.fig}.  Then the corresponding toric variety is the
projective plane $\Proj ^{2}$.  Assume that $\deg $ takes the
indicated values on the generators of the one-dimensional cones.  Then
after simplifying, one sees that the toric form associated to this
data is
\[
f_{N,\deg} (q) = \sum _{a,b\in \ZZ }\frac{1-\ee^{2\pi\ii\alpha
}\ee^{2\pi\ii \beta }\ee^{2\pi\ii \gamma }}{(1-\ee^{2\pi\ii\alpha }q^{a})
(1-\ee^{2\pi\ii\beta }q^{b})(1-\ee^{2\pi\ii\gamma }q^{-a-b})}.
\]
\begin{figure}[ht]
\psfrag{a}{$\alpha $}
\psfrag{b}{$\beta $}
\psfrag{c}{$\gamma $}
\centerline{\includegraphics[scale = .3]{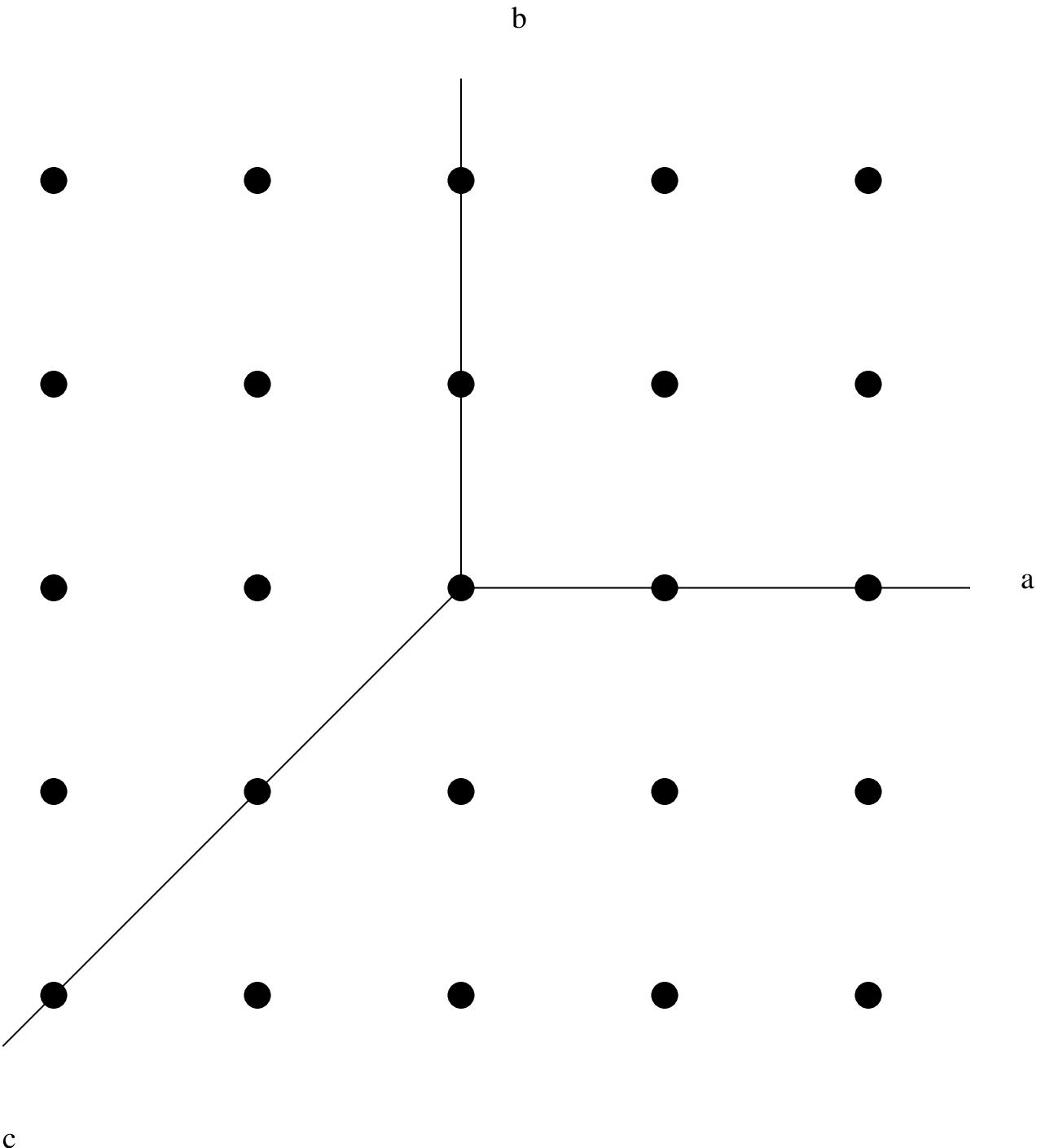}}
\caption{\label{p2.fig}A toric form associated to $\Proj ^{2}$.}
\end{figure}

\end{example}

In the rest of this section we will show that $f_{N,\deg}$ is well-defined.
Moreover, in \S\ref{toric.ell}, we will justify our nomenclature by
showing that if $\deg (N)\subset \QQ$, then $f_{N,\deg}$ is a
holomorphic modular form of weight $r$.

\subsection{}
Our first goal is to investigate the function 
$$
r(q,m)=\sum_{C\in\Sigma}(-1)^{\codim C} 
\ac{\sum_{n\in C} q^{m\cdot n} \ee^{2\pi\ii\deg(n)}}.
$$

First of all, it is easy to see that $r(q,m)$ does not depend on
$\Sigma $, in the sense that $\Sigma $ can be replaced by any fan on
which $\deg $ is piecewise-linear. Indeed, it suffices to show that
$r(q,m)$ does not change if $\Sigma $ is subdivided.  This is easily
seen to be true for those cones and $m\in M$ where the corresponding
sums are absolutely convergent.  Analytic continuation then finishes
the argument.  Therefore, $f_{N,\deg}$ is independent of $\Sigma$, as
our notation suggests.

\subsection{}
Now we want to give a homological interpretation of $r(q,m)$.  Let
$\CC[N]^\Sigma$ be the following deformation of $\CC[N]$. As a vector
space $\CC[N]^\Sigma = \otimes_{n\in N} \CC y^n$, where $y$ is a dummy
multi-variable.  We define multiplication in $\CC [N]^{\Sigma }$ by
$$
y^{n_1}y^{n_2}=
\left\{
\begin{array}{ll}
y^{n_1+n_2} & \hbox{if there exists $C\in\Sigma$ with $n_1,n_2\in C$,}\\
0 & \hbox{otherwise.}
\end{array}
\right.
$$
Furthermore, for each cone $C\in\Sigma$, we have the submodule $\CC[C]
\subset \CC [N]^{\Sigma }$.

Let $R$ be the polynomial ring $\CC[y^{d}\mid d \in D]$, where $D$ is the set
of generators of the one-dimensional cones of $\Sigma$. There is a
natural action of $R$ on $\CC[N]^\Sigma$, where any product
$y^{d}y^{n}$ is zero unless $d$ and $n$ lie in some cone of
$\Sigma$. It is straightforward to see that $\CC[N]^\Sigma$ and $\CC
[C]$ are finitely generated $R$-modules.

The modules $\CC [N]^{\Sigma }$ and $\CC [C]$ have a natural grading
by $N$ that is compatible with the $R$-action.  Moreover, these
modules have a $\CC $-grading induced by $\deg$.  In particular, the
spaces 9$\Tor^i_R(\CC[N]^\Sigma,\CC)$ have a natural $N\oplus\CC$
grading.

\begin{proposition}
Let $D$ be the set of generators of one-dimensional cones of $\Sigma $.
Then
\begin{equation*}
r(q,m)=
{\sum_i (-1)^i {\Trace}_{\Tor^i_R(\CC[N]^\Sigma,\CC)}
q^{m\cdot {\bf n}} \ee^{2\pi\ii {\bf deg}}
\over
\prod_{d\in D}(1-\ee^{2\pi\ii \deg(d)}q^{m\cdot d})
}.
\end{equation*}
Here $m\cdot{\bf n}$ (respectively $\bf deg$) denotes the linear
operator that takes the 
value $m\cdot n$ (resp. $\deg$) on each graded component.
\label{syzygies}
\end{proposition}

\begin{proof}
To begin, notice the long exact sequence of $\Tor $ implies that the
function
\begin{equation}\label{rsubA}
r(A;q,m) := 
{\sum_i (-1)^i {\Trace}_{\Tor^i_R(A,\CC)}
q^{m\cdot {\bf n}} \ee^{2\pi\ii {\bf deg}}
\over
\prod_{d\in D}(1-\ee^{2\pi\ii \deg(d)}q^{m\cdot d})
}
\end{equation}
is additive on graded $R$-modules $A$. For $i=0,\dots ,r$, let
$C^{i}\subset \Sigma $ be the set of cones of codimension $i$.
Consider the complex of $R$-modules
$$
0\longrightarrow\bigoplus_{C\in C^{0}}\CC[C]
\longrightarrow\bigoplus_{C\in C^{1}}\CC[C]
\longrightarrow\cdots
\longrightarrow\bigoplus_{C\in C^{r}}\CC[C]\longrightarrow 0,
$$
where the differentials are induced by the boundary maps from any cone $C$ to
its maximum proper subcones, with appropriate signs depending on the
relative orientations of $C$ and these subcones. It is straightforward
to see that the cohomology of this complex occurs only at 
the $C^{0}$ spot, and in fact equals $\CC[N]^\Sigma$.
Since $r(A;q,m)$ is additive, this means that we only need to show 
$$
r(\CC[C];q,m) = \ac{\sum_{n\in C} q^{m\cdot n} \ee^{2\pi\ii \deg(n)}}
$$
for each cone $C\in\Sigma$.

So fix $C\in \Sigma $.  Let $D (C)$ be the set of generators of
one-dimensional cones of $C$, and let $R_{C}$ be the polynomial ring
$\CC[y^{d}\mid d\in D (C)]$.  We claim that in the denominator of
(\ref{rsubA}), we need to take the product only over $d\in D (C)$.
Indeed, if $d\not \in D (C)$, then $y^{d}$ acts trivially on $\CC[C]$.
Hence $\Tor^*_{R}(\CC[C],\CC)$ is equal to $\Tor^*_{R_{C}}
(\CC[C],\CC)$ tensored with the Koszul complex for $\{y^{d}\mid
d\notin D (C)\}$.  This Koszul complex gives rise to extra factors in
the alternating sum of traces, which will cancel the factors
\[
(1-\ee^{2\pi\ii \deg(d)}q^{m\cdot d}),\quad d\not \in D (C)
\]
in the denominator.

Hence we must now show that
\begin{equation}\label{final.statement}
{
\sum_i (-1)^i {\Trace}_{\Tor^i_{R_C}(\CC[C],\CC)}
q^{m\cdot {\bf n}} \ee^{2\pi\ii {\bf deg}}
\over
\prod_{d\in D (C)}(1-\ee^{2\pi\ii \deg(d)}q^{m\cdot d})
}
= \ac{\sum_{n\in C}
q^{m\cdot n}
\ee^{2\pi\ii \deg(n)}}.
\end{equation}
Let
$m\in M$ satisfy (\ref{cond.on.m}), so that the series on the right of
\eqref{final.statement} converges absolutely for all $|q|<\varepsilon$
for some $\varepsilon>0$. We will prove the more general statement
\begin{equation}\label{iden}
{
\sum_i (-1)^i {\Trace}_{\Tor^i_{R_C}(A,\CC)}
q^{m\cdot {\bf n}} \ee^{2\pi\ii {\bf deg}}
\over
\prod_{d\in D (C)}(1-\ee^{2\pi\ii \deg(d)}q^{m\cdot d})
}
= {\Trace}_{A}q^{m\cdot {\bf n}}\ee^{2\pi\ii {\bf deg}}
\end{equation}
for all finitely generated graded $R_C$-modules $A$.  Indeed, in
\eqref{iden} both sides are additive on finitely generated graded
$R_C$-modules $A$, and coincide on $R_C$ itself.  Together with the
existence of a free resolution, this implies \eqref{iden}, which
finishes the proof of the lemma.
\end{proof}

\subsection{}
The following important lemma will allow us to show convergence
of the series used to define $f_{N,\deg}$. 
\begin{lemma}
Let $n\in N$ be such that $q^{m\cdot n}$ appears with nonzero
coefficient in
\begin{equation}\label{numer}
\sum_i (-1)^i {\Trace}_{\Tor^i_R(\CC[N]^\Sigma,\CC)}
q^{m\cdot {\bf n}} \ee^{2\pi\ii {\bf deg}}.
\end{equation}
Let $D = \{d_{i} \}$ 
be set of generators of one-dimensional cones of $\Sigma $.
Then $n$
lies in the interior of the convex hull of the finite set of points
$$
\Bigl\{\sum_{i\in I}d_i \Bigm | \hbox{$I$ is any subset of the set
of all $i$}\Bigr\}.
$$
\label{convexhull}
\end{lemma}

\begin{proof}
Without loss of generality, we may assume that $\Sigma$ is simplicial.
Indeed, we can freely subdivide $\Sigma$ as long as we do not add any
new one-dimensional cones. If we pick a generic collection of points
on the one-dimensional faces of $\Sigma$, and construct the convex
hull of these points for each cone, then we construct a simplicial
refinement of $\Sigma$.

It is more convenient now to equip $\CC[N]^\Sigma$ with a new
$\QQ^{\sharp D}_{\geq 0}$-grading that we will denote $\widehat\deg$.
To define this grading, suppose that $n\in N$ lies in the cone
$C\in\Sigma$, and write 
\[
n = \sum _{d\in D (C)} \alpha _{d}\cdot d, \quad \alpha _{d}\in \QQ.
\]
Then $\widehat\deg ( n)$ is defined to have component $\alpha _{d}$ for
$d\in D (C)$, and component $0$ otherwise.  Obviously, this definition
is independent of $C$, and it is easy to see that the old $N\oplus
\CC$ grading is induced from this finer grading.  Moreover,
$\CC[N]^\Sigma$ is a graded $R$-module in the new sense for
$R=\CC[y^{d}\mid d\in D]$.

Suppose that $n$ is a lattice point that violates the statement of the
lemma. This means that $n$ appears with a nonzero coefficient in
\eqref{numer}, and that there exists $h\in M_\CC$ such that 
\begin{equation}\label{cond}
h\cdot n \geq h\cdot(\sum_{i\in I}d_i),\quad \hbox{for all $I \subset
\{1,\dots ,\sharp D\}$.}
\end{equation}
Using the notation $n_I:=n-\sum_{i\in I} d_i$, the 
condition \eqref{cond} translates to 
\begin{equation}\label{cond2}
h\cdot n_I\geq 0, \quad\hbox{for all $I$.} 
\end{equation}

To calculate $\Tor^*(\CC[N]^\Sigma,\CC)$, we can tensor
$\CC[N]^\Sigma$ with the Koszul complex and calculate the cohomology of
the resulting complex.  At the $N$-degree $n$ we will have
\begin{equation}\label{koz.cpx}
0\longrightarrow\cdots\longrightarrow
\bigoplus_{i<j}(\CC[N]^\Sigma[-d_i-d_j])_n\longrightarrow
\bigoplus_{i}(\CC[N]^\Sigma[-d_i])_n\longrightarrow
(\CC[N]^\Sigma)_n\longrightarrow 0.
\end{equation}
For each $I$ the space $(\CC[N]^\Sigma[-\sum_{i\in I}d_i])_n$
is one-dimensional, with $\widehat\deg$ equal to
$$
\widehat\deg(n_I)+\sum_{i\in I}\widehat\deg(d_i).
$$

Now the complex \eqref{koz.cpx} splits into a direct sum of
subcomplexes according to $\widehat\deg$.  If $n$ appears with a
nonzero coefficient in \eqref{numer}, then the Euler characteristic of one
of these subcomplexes, say $\CCC _{0}$, must be nonzero.  We will
show that if this $n$ satisfies \eqref{cond2}, then in fact $\chi (\CCC
_{0}) = 0$, which is a contradiction.  

Let $\widehat\deg_0$ be the degree of $\CCC _{0}$.  Consider the set
$\III$ of all $I$ such that
$$
\widehat\deg(n_I)+\sum_{i\in I}\widehat\deg(d_i) = \widehat\deg_0.
$$  
We claim that $\III$ has a unique maximal element with respect to inclusion.
Indeed, suppose $I, J \in \III$, so that 
$$
\widehat\deg(n_I)+\sum_{i\in I}\widehat\deg(d_i) = 
\widehat\deg(n_J)+\sum_{j\in J}\widehat\deg(d_j).
$$
Then
$$
\widehat\deg(n_I)-\widehat\deg(n_J)=
\sum_{j\in (J\smallsetminus(J\cap I))}\widehat\deg(d_j)-
\sum_{i\in (I\smallsetminus(I\cap J))}\widehat\deg(d_i),
$$
which implies that $d_j$ lies in the minimum cone containing $n_I$
for all $j\in (J\smallsetminus(J\cap I))$, and that the corresponding
coefficient of 
$n_{I}$ is $\geq 1$.  Hence
$$
\widehat\deg(n_{I\cup J})+\sum_{j\in (J\smallsetminus(J\cap
I))}\widehat\deg(d_i) = \widehat\deg(n_I),
$$
and
$$
\widehat\deg(n_I)+\sum_{i\in I}\widehat\deg(d_i) = 
\widehat\deg(n_{I\cup J})+\sum_{j\in (I\cup J)}\widehat\deg(d_j).
$$
Hence $I\cup J \in \III$.
We will denote the maximum set for $\widehat\deg_0$ by 
$I^{\max}$. 

It is easy to describe the subsets $I \subset I^{\max}$
satisfying
\begin{equation}\label{allowed}
\widehat\deg(n_{I^{\max}\smallsetminus I})+
\sum_{i\in (I^{\max}\smallsetminus I)}\widehat\deg(d_i) =
\widehat\deg_0.
\end{equation}
Indeed, \eqref{allowed} happens if and only if there exists a cone
$C\in\Sigma$ containing both $n_{I^{\max}}$ and $\{d_i\mid i\in I\}$.
Let us say that $I$ is an \emph{allowed subset} of $I^{\max}$ if this
happens.

Let $C_{\min}$ be the minimal cone containing $n_{I^{\max}}$.  We will
now show that if $i\in I^{\max}$, then $d_{i}$ is not a generator of
$C_{\min}$.  Indeed, suppose $d_{i}$ is a generator of $C_{\min}$, and
let $I$ be an allowed subset.  If $i\not \in I$, then $I\cup \{i \}$
is allowed, and if $i\in I$, then $I\smallsetminus \{i \}$ is allowed.
Hence the set of all allowed subsets splits into pairs, and each pair
contributes $0$ to $\chi (\CCC_0)$.  However, by assumption $\chi
(\CCC _{0})\not =0$.

In addition, if any of the coordinates of $\widehat\deg(n_{I^{\max}})$
were at least $1$, then we could replace $I^{\max}$ with $I^{\max}\cup
\{i \}$, which contradicts its maximality.  This means that
$n_{I^{\max}}$ lies in the interior of the convex hull of all possible
$\sum_{j\in J} d_j$, where $J$ satisfies
\[
J\subset \{i\mid \hbox{$d_{i}\in D ( C_{\min})$} \}.
\]
Therefore, $0$ is in the interior of the convex hull of
all $n_{I^{\max}\cup J}$, which implies
$$
h\cdot n_{I^{\max}} = 0\quad \hbox{and}\quad h\cdot d_j=0
$$
for all $d_j\in D (C_{\min})$.  Because $h\cdot n_I\geq 0$ for all
$I$, we conclude that $h\cdot d_i$ is nonnegative for all $i\in
I^{\max}$, and is nonpositive for all other $i$.

Let us now calculate the Euler characteristic of ${\CCC}_0$.  We have
\begin{align}
(-1)^{\sharp(I^{\max})} \chi({\CCC}_0) &= 
\sum_{I\subseteq I^{\max}, I~{\rm allowed}}(-1)^{\sharp(I)}\\
\label{two}&=
\sum_{C\in\Sigma,C\ni n_{I^{\max}}}(-1)^{\codim C}
\sum_{I\subseteq  I^{\max},d_i\in C ~{\rm for~all}~i\in I}
(-1)^{\sharp(I)}.
\end{align}
Here we have used $\sum_{C\in\Sigma,C\ni n}(-1)^{\codim C} = 1$ for
every lattice point $n$.  This is easy to verify by looking at the
simplicial complex $\Sigma $ induces on a small sphere about $n$.

We can quotient $N_{\RR}$ by the subspace generated by $C_{\min}$,
and can induce a fan there.  Hence without loss of generality we may
assume $n_{I^{\max}} = 0$.

We now observe that the inner sum in \eqref{two} is zero, unless no $d_i$ with
$i\in I^{\max}$ lies in $C$.  To see what this means, define 
\[
\tilde{\Sigma } = \left\{C\in \Sigma \bigm| \hbox{there exists
$d_{i}\in D (C)$ with $i\in I^{\max}$} \right\}.
\]
Then
\begin{align}
(-1)^{\sharp(I^{\max})} \chi({\CCC}_0) &= 
\sum_{\substack{C\not \in \tilde{\Sigma}}}
(-1)^{\codim C}\\
&= 1 - \sum_{\substack{C\in\tilde{\Sigma }}} 
(-1)^{\codim C}\label{2nd}.
\end{align}

Now define a subset $U$ by 
\[
U := \bigcup _{C\in \tilde{\Sigma }} \relint (C).
\]
We claim that the sum in \eqref{2nd}
is the Euler characteristic of $U$ in \v{C}ech cohomology theory.
Indeed, this cohomology could be calculated using the acyclic covering
by the open stars 
\[
{\rm Star}(d_i):=\bigcup_{C_1\in\Sigma, C_1\ni d_i}\relint (C_{1}),
\]
for $i\in I^{\max}$, 
and the cones $C$ above are in one-to-one correspondence with
intersections of ${\rm Star}(d_i)$ via 
\[
{\rm Star}(C):=\bigcup_{C_1\in\Sigma, C_1\supseteq C}\relint (C_{1}).
\]
Therefore, to show that $\chi({\CCC}_0)=0$, it suffices to
show that the space $U$ is contractible. 

First of all, notice that $U$ contains the entire open half-space
$h>0$. Pick a point $p$ inside this half-space. For every point $u\in
U$ we will define a path $\nu\colon [0,1]\to U$ such that $\nu ( 0)
=u$, $\nu (1)= p$ as follows.  Consider any cone $C$ with $u\in
C$. There is a unique decomposition
$$u=u_1+u_2$$
such that $u_1$ has nonzero coordinates only for $d_i\in I^{\max}$
and $u_2$ has nonzero coordinates only for $d_i\notin I^{\max}$.
Then define 
\[
\nu (t) = 
\begin{cases}
u_1 + (1-2t)u_2,&t\in [0,1/2],\\
(2-2t)u_1  +  (2t-1)p,&t\in [1/2,1].
\end{cases}
\]
It is easy to see 
that these paths assemble together into a continuous map
$U\times[0,1]\to U$ providing a retraction of $U$ onto $p$.

Hence we have shown that $\chi({\CCC}_0)=0$, which contradicts our
assumption that $n$ and $\widehat\deg_0$ disobey the conditions of 
the lemma.  This completes the proof.
\end{proof}

\begin{theorem}
The function $f_{N,\deg}(q)$ of Definition \ref{modularfromdegree}
is well-defined.  More precisely,
\begin{enumerate}
\item  There exists $\varepsilon>0$ such that the outer sum over $M$
converges absolutely and uniformly for all $|q|<\varepsilon$.
\item  Each term $r(q,m)$ has a $q$-expansion with nonnegative powers
of $q$, and only a finite number of $m$ contribute nontrivially to
$q^n$ for any fixed $n$.
\end{enumerate}
\label{welldefd}
\end{theorem}

\begin{proof}
The lattice $M$ is a disjoint union of a finite number of 
regions, each characterized by the collection of signs of
$m\cdot d_i$ (positive, zero, negative) for all $i$. Each such region
is the interior of a rational polyhedral cone in some sublattice of 
$M$. It therefore suffices to show that each of the above convergence
statements is true separately for $m$ inside one such region $\RRR$. 

Pick $\varepsilon<1$ such that 
\[
\varepsilon < |\ee^{2\pi\ii\deg(d_i)}|^{-1}/2, \quad d_{i}\in D.
\]
Then for every $d_i$ with ${\RRR}\cdot d_i>0$, 
we have 
$$
|1-\ee^{2\pi\ii\deg(d_i)}q^{m\cdot d_i}|>1/2.
$$
For any $d_i$ with ${\RRR}\cdot d_i=0$ we just get a nonzero constant
in the denominator, because we always assume that $\deg(d_i)$ are
not integers. Finally, Lemma \ref{convexhull} implies that for every generator $m_j$ of the closure of 
$\RRR$, and every $n$ contributing to the numerator of $r(q,m)$
we get 
$$m_i\cdot n > \sum_{k,m_i\cdot d_k<0} m_i\cdot d_k.$$
Since for $m\cdot d_i<0$ and $|q|<\varepsilon$, we have 
$$|1-\ee^{2\pi\ii\deg(d_i)}q^{m\cdot d_i}|>c_1q^{m\cdot d_i}$$
for each $m$, it follows that the terms of the $\sum_m r(q,m)$ can be estimated
by $c_2 q^{l(m)}$,  where $l(m)$ is the minimum of the linear functions
(one for each $n$ in the numerator of $r(q,m)$) that 
are strictly positive on all generators of the closure of $\RRR$.
This implies the first statement.  The second statement is treated
similarly.
\end{proof}

%
%

\section{Theta functions and toric forms}\label{chern.sect}

\subsection{}
The space $\TTT $ has a natural ring structure,
with multiplication given by  
$$
f_{N_1,\deg_1}
f_{N_2,\deg_2}
=
f_{N_1\oplus N_2,\deg_1\oplus \deg_2}.
$$

In this section we study the ring $\TTT $ by expressing the toric form
$f_{N,\deg }$ in terms of the theta function with characteristic
\cite{Chandra}
\begin{equation}\label{theta.def}
\theta _{11} (z,\tau ) = \theta (z,\tau ) :={1\over
i}\sum_{n\in\ZZ}(-1)^n \ee^{\pi\ii\tau(n+{1\over 2})^2} \ee^{\pi \ii z
(2n+1)}.
\end{equation}

To do this, we give a topological interpretation of $f_{N,\deg }$.
Associated to the pair $(N,\deg )$ is complete rational polyhedral fan
$\Sigma \subset N\otimes \RR$, and hence a complete toric variety $X$.
Suppose that $X$ is nonsingular.  Then we show in Theorem~\ref{klem}
that $f_{N,\deg }$ can be computed as the graded Euler characteristic
of certain graded infinite-dimensional vector bundle $W$ over $X$.
Moreover, we can compute this Euler characteristic $\chi (W)$ in two
different ways:
\begin{enumerate}
\item We can use the \v{C}ech cohomology complex associated to the
covering of $X$ by toric affine charts.  This yields the
usual expression for $f_{N,\deg }$.
\item We can use the Hirzebruch-Riemann-Roch Theorem to compute $\chi
(W)$ in terms of the Todd class $\Td (X)$ and the Chern character
$\chernchar (W)$.  Using the Jacobi triple product formula 
\begin{equation}\label{jf}
\theta(z,\tau)=q^{1 \over 8}  (2 \sin \pi z)
\prod_{l=1}^{l=\infty}(1-q^l)
 \prod_{l=1}^{l=\infty}(1-q^l \ee^{2 \pi \ii z})(1-q^l \ee^{-2 \pi \ii
z}),
\end{equation}
we can interpret the resulting expressions in terms of $\theta $.
\end{enumerate}
This technique was used by Borisov-Libgober \cite{borlibg} to
compute the elliptic genus of a smooth toric variety, and more details
can be found there. 

More generally, if $X$ is singular, we can compute these quantities
through a limiting process (Theorem~\ref{limit}).  As a corollary of
these results, we obtain that $f_{N,\deg}$ is meromorphic in $\HHH $.

\subsection{}
We begin with some notation.  Let $\deg $ be a degree function with
respect to the complete fan $\Sigma $, and let $\{d_{i} \}$ be the set
of generators of one-dimensional cones of $\Sigma $.  Let $X$ be the
toric variety associated to $\Sigma $, and for each $d_{i}$, let
$D_{i}\subset X$ be the corresponding toric divisor.  In what
follows, we will usually abuse notation and use $D_{i}$ to mean either
the divisor or its cohomology class.

Let $\OOO (D_{i})$ be the line bundle associated to $D_{i}$.  Recall
that the space of sections of $\OOO(D_{i})$ over each open subset of
$X$ is by definition the space of all rational functions on $X$ that
have a pole of order at most $1$ along $D_{i}$ and no other poles.
Moreover, we may extend the torus action from $X$ to each $\OOO
(D_{i})$ by inducing from the action on the field of rational
functions.

\subsection{}

Let $\OOO $ be the trivial line bundle.  For any $\alpha \in \CC$, we
define graded bundles 
\begin{align*}
\Lambda_{\alpha }^{*} (\OOO (D_{i})) &:= \OOO \oplus \OOO _{(\alpha )} (D_{i})\\
\Sym_{\alpha }^{*} (\OOO (D_{i})) &:= \OOO \oplus \OOO _{(\alpha )}
(D_{i}) \oplus \OOO _{(\alpha^{2} )} (2D_{i}) \oplus \cdots
\end{align*}

For each generator $d_{i}$ let $\alpha _{i} = \deg (d_{i})$.
Finally let $W$ be the graded bundle 
\[
\bigotimes _{i}\bigotimes _{n=1}^{\infty }\left (
\Lambda^*_{-q^n\ee^{-2\pi\ii\alpha_i}}{\OOO}(D_i) \otimes
\Lambda^*_{-q^{n-1}\ee^{2\pi\ii\alpha_i}}{\OOO}(-D_i) \otimes
{\Sym}^*_{q^{n}}{\OOO}(D_i) \otimes
{\Sym}^*_{q^{n}}{\OOO}(-D_i)\right).
\]
Although the bundle $W$ is infinite-dimensional, each graded piece
corresponding to a fixed power of $q$ is finite-dimensional.

We are now ready to state the main result of this section.

\begin{theorem}
Assume that the toric variety $X$ is nonsingular, and that $\alpha
_{i}\not \in \ZZ$ for all generators of one-dimensional cones of
$\Sigma $.  Then
\begin{equation}\label{state}
f_{N,\deg}(q)=\int_X \prod_i {(D_i/2\pi\ii)\theta(D_i/2\pi\ii-\alpha
_{i},\tau )\theta'(0,\tau ) \over \theta(D_i/2\pi\ii,\tau
)\theta(-\alpha_{i},\tau ) }.
\end{equation}
Here the derivative is the partial derivative with respect to $z$, the
product is taken over all $d_{i}$, and the right hand side of
\eqref{state} is interpreted as an expression in $H^{*} (X)$ using the
Jacobi triple product \eqref{jf}.  \label{klem}
\end{theorem}

\begin{proof}
First consider the right-hand side of \eqref{state}.  Using
$${\Td}(X) = \prod_i{D_i\over 1- \ee^{-D_i}},$$
and applying \eqref{jf}, and the product formula for $\theta' (0,\tau
) $ \cite{Chandra},
the right-hand side becomes
\[
\int_X \prod_i {D_i\over (1-\ee^{-D_i})}\cdot \prod_{n=1}^{\infty}
{(1-q^n)^2(1-q^n\ee^{D_i-2\pi\ii\alpha_i})
(1-q^{n-1}\ee^{-D_i+2\pi\ii\alpha_i}) \over
(1-q^n\ee^{D_i})(1-q^n\ee^{-D_i})
(1-q^n\ee^{-2\pi\ii\alpha_i})(1-q^{n-1}\ee^{2\pi\ii\alpha_i}) }.
\]

This is seen to be
$$
\prod_i\prod_{n=1}^\infty
{
(1-q^n)^2
\over
(1-q^n\ee^{-2\pi\ii\alpha_{i}})
(1-q^{n-1}\ee^{2\pi\ii\alpha_{i}})
}
\cdot
\int_X \Td(X) 
{\chernchar}
(W),
$$
which by Hirzebruch-Riemann-Roch is
$$
=\prod_i\prod_{n=1}^\infty
{
(1-q^n)^2
\over
(1-q^n\ee^{-2\pi\ii\alpha_{i}})
(1-q^{n-1}\ee^{2\pi\ii\alpha_{i}})
}
\cdot
\chi(W).
$$

Now we want to calculate $\chi (W)$ using the \v{C}ech complex
associated to the covering of $X$ by toric affine charts.  Let $x$ be
a dummy multi-variable that keeps track of the torus action on $\OOO
(D_{i})$.

For every cone $C\in \Sigma$ with generators $d_1,\dots ,d_k$, let
$U_{C}$ be the corresponding affine chart.  Let 
\[
\Dim_{M}(\Gamma(W,U_C)) = \sum _{m\in M}\Dim \Gamma (
W,U_C)_{m}\cdot x^{m}
\]
be the graded dimension of $\Gamma(W,U_C)$ with respect to $M$.  Let
$m_1,\dots ,m_k$ be the dual basis to $d_1,\dots ,d_k$.  Then the
space of sections of ${\OOO}(D_j)$ on $U_C$, as a module over
$\Gamma(\OOO ,U_C)$, is generated by an element with grading $-m_j$ if $d_j$
is one of the generators of $C$, and by an element with zero grading
otherwise. Now it is easy to see that
$$
\prod_i\prod_{n=1}^\infty
{
(1-q^n)^2
\over
(1-q^n\ee^{-2\pi\ii\alpha_{i}})
(1-q^{n-1}\ee^{2\pi\ii\alpha_{i}})
}
\cdot
\Dim_{M}(\Gamma(W,U_C))
$$
$$=
(\sum_{m\in M,m\cdot C=0}x^m)\cdot
\prod_{i,d_i\in C}\prod_{n=1}^\infty
{
(1-q^n)^2
(1-q^n\ee^{-2\pi\ii\alpha_i} x^{-m_i})
(1-q^{n-1}\ee^{2\pi\ii\alpha_i} x^{m_i})
\over
(1-q^n\ee^{-2\pi\ii\alpha_i})(1-q^{n-1}\ee^{2\pi\ii\alpha_i})
(1-q^n x^{-m_i})(1-q^{n-1}x^{m_i})
}
$$
$$
=(\sum_{m\in M,m\cdot C=0}x^m)\cdot
\prod_{i,d_i\in C}
(
\sum_{l\in \ZZ}{x^{lm_i}\over (1-\ee^{2\pi\ii\alpha_i}q^l)}
)
=
\sum_{m\in M}x^m\prod_{i,d_i\in C} {1\over
(1-\ee^{2\pi\ii\alpha_i}q^{m\cdot d_i})},
$$
where we expand everything as a power series around $q=0$.
Here we have used the identity of power series in $q$
$$
\prod_{k\geq 1}
{(1-tyq^{k-1})(1-t^{-1}y^{-1}q^k)
\over
(1-tq^{k-1})(1-t^{-1}q^k)
}
=
\sum_{m\in \ZZ} \left({t^m\over 1-yq^m}\right)
\prod_{k\geq 1}  
{(1-yq^{k-1})(1-y^{-1}q^k)\over(1-q^k)^2},
$$
whose proof can be found in \cite{borlibg}.

Now all we need to prove the statement is to notice that, as power
series in $q$,
$$
\chi(W)=\sum_{C\in \Sigma}(-1)^{\codim C}\Dim_{M}(\Gamma(W,U_C))\bigr|_{x=1}.
$$
This follows easily from the description of \v{C}ech complex,
see \cite{borlibg} for details. 
\end{proof}

We will now combine the results of Theorems~\ref{welldefd}
and~\ref{klem} to get similar expressions for an arbitrary toric form.

\begin{theorem}
Let $\deg$ be a degree function with respect to a complete simplicial
fan $\Sigma$ in a lattice $N$.  Let $X$ be the toric variety
associated to $\Sigma $, and $\widehat\Sigma$ be a refinement of
$\Sigma$ providing a desingularization $\widehat X$ of $X$. Let
$\deg_1$ be a generic degree function with respect to $\Sigma$. Then
\begin{equation*}
f_{N,\deg}(q)=
\lim_{\varepsilon\to 0}
\int_{{\widehat X}} \prod_i
{(D_i/2\pi\ii)\theta(D_i/2\pi\ii-\deg(d_i)-\varepsilon \deg_1(d_i),\tau )\theta'(0,\tau )
\over
\theta(D_i/2\pi\ii,\tau )\theta(-\deg(d_i)-\varepsilon \deg_1(d_i),\tau )
},
\end{equation*}
where the product is taken over all generators of one-dimensional cones of
$\widehat\Sigma$.
\label{limit}
\end{theorem}

\begin{proof}
The point is that sometimes when we desingularize and add extra $d_i$,
we may get $\deg(d_i)\in\ZZ$, so Theorem~\ref{klem} cannot be directly
applied. On the other hand, it is easy to see that in Theorem
\ref{welldefd}, the convergence is uniform in $\deg$.  Hence
$f_{N,\deg}$ is continuous as a function of $\deg$, and the result
follows.
\end{proof}

\begin{corollary}
The function $f_{N,\deg}$ is meromorphic on the upper
halfplane. \qed
\end{corollary}

%
%

\section{Toric forms of level $l$}\label{toric.ell}

\subsection{}
Let $l\geq 2$ be an integer, and let $\MMM _{r}(\Gamma_{1} (l) )$ be
the ring of modular forms of weight $r$ with respect to the congruence
subgroup $\Gamma _{1} (l) \subset SL_{2} (\ZZ )$.  Let $\TTT
_{r}(l)\subset \TTT_{r} $ be the subspace spanned by
\[
\Bigl\{f_{N,\deg }\in \TTT _{r}\Bigm | \deg(N)\subset {1\over l}\ZZ\Bigr\}.
\]
We say that $\TTT (l) = \bigoplus _{r} \TTT _{r} (l)$ is the ring of
toric forms of \emph{level $l$}.

In this section, we show that $\TTT _{r}(l)\subseteq \MMM
_{r}(\Gamma_{1} (l) )$, and describe an explicit set of generators for
$\TTT (l)$.  We relate toric forms of level $l$ to the elliptic genera
considered by Hirzebruch.

\subsection{}
To begin, we show that toric forms of level $l$ are automorphic with
respect to $\Gamma _{1} (l)$.  For later purposes, it is more
convenient to phrase the statement in terms of $\Gamma _{0} (l)$,
the subgroup of matrices satisfying
\[
\Bigl (\begin{array}{cc}
a&b\\
c&d
\end{array}\Bigr)=
\Bigl (\begin{array}{cc}
*&*\\
0&*
\end{array}\Bigr)\mod l.
\]

\begin{proposition}\label{modularity.of.tf}
If $f_{N,\deg} \in \TTT_{r} (l)$ and 
$\left ( \begin{array}{cc}a&b\\
c&d\end{array}\right)\in \Gamma _{0} (l)$, then
\[
f_{N,\deg}\left (\frac{a\tau +b}{c\tau +d}\right) = (c\tau
+d)^{r}f_{N,d\cdot \deg} (\tau ).
\]
\end{proposition}

\begin{proof}
This will follow from Theorems~\ref{klem} and~\ref{limit} and the
transformation properties of $\theta $.
For any modular transformation $\left ( \begin{array}{cc}a&b\\
c&d\end{array}\right)\in SL_{2} (\ZZ )$, we have
\begin{align*}
\theta\Bigl({z\over 
c\tau+d},{a\tau+b\over c\tau+d}\Bigr)&=\zeta (c\tau+d)^{1\over2}
\ee^{\pi \ii cz^2\over c\tau+d}\theta(z,\tau)\\
\theta'\Bigl(0,{a\tau+b\over c\tau+d}\Bigr)&=\zeta (c\tau+d)^{3\over2}
\theta'(0,\tau),
\end{align*}
where $\zeta ^{8} =1$, and $\zeta $ depends on $a,b,c,d$ but not on
$\tau$ or $z$ \cite{tanner}.  
Now use the notation of Theorem~\ref{limit}.  We denote $\deg (d_{i})$
by $\alpha _{i}$ and $\deg_{1} (d_{i})$ by $\beta _{i}$.  We have 
\begin{align*}
&f_{N,\deg} \left(\frac{a\tau +b}{c\tau +d} \right) (c\tau +d)^{-r} =
\lim _{\varepsilon \rightarrow 0}\int _{\hat{X}} \prod
_{i}\frac{(\frac{D_{i}/2\pi \ii }{c\tau +d})\theta (\frac{D_{i}/2\pi
\ii }{c\tau +d} - \alpha _{i}-\varepsilon \beta _{i}, \frac{a\tau +b}{c\tau +d} )\theta ' (0,\frac{a\tau +b}{c\tau +d} )}{\theta (\frac{D_{i}/2\pi
\ii }{c\tau +d}, \frac{a\tau +b}{c\tau +d})\theta (-\alpha _{i}-\varepsilon \beta _{i},\frac{a\tau +b}{c\tau +d})} \\
&=\lim _{\varepsilon \rightarrow 0}\int _{\hat{X}} \prod
_{i}\frac{(\frac{D_{i}}{2\pi \ii})\ee ^{\frac{\pi \ii c }{c\tau +d}(\frac{D_{i}}{2\pi \ii }- (\alpha
_{i}+\varepsilon \beta _{i}) (c\tau +d))^{2}} \theta (\frac{D_{i}}{2\pi \ii}- (\alpha
_{i}+\varepsilon \beta _{i}) (c\tau +d), \tau) \theta '(0,\tau
)}{\theta (\frac{D_{i}}{2\pi \ii},\tau )\theta (-(\alpha _{i}+\varepsilon
\beta _{i}) (c\tau +d),\tau )\ee ^{\frac{\pi \ii c}{c\tau +d} ((\frac{D_{i}}{2\pi \ii })^{2} + (\alpha
_{i}+\varepsilon \beta _{i})^{2} (c\tau +d)^{2})}}\\
&=\lim _{\varepsilon \rightarrow 0}\int _{\hat{X}} \prod
_{i}\frac{(\frac{D_{i}}{2\pi \ii})\theta (\frac{D_{i}}{2\pi \ii}- (\alpha
_{i}+\varepsilon \beta _{i}) (c\tau +d), \tau) \theta '(0,\tau
)}{\theta (\frac{D_{i}}{2\pi \ii},\tau )\theta (-(\alpha _{i}+\varepsilon
\beta _{i}) (c\tau +d),\tau )}\ee ^{-c D_{i} (\alpha _{i}+\varepsilon \beta _{i}) }\\
&=\lim _{\varepsilon \rightarrow 0}\int _{\hat{X}} \prod
_{i}\frac{(\frac{D_{i}}{2\pi \ii})\theta (\frac{D_{i}}{2\pi \ii}- d\alpha
_{i}-\varepsilon \beta _{i} (c\tau +d), \tau) \theta '(0,\tau
)}{\theta (\frac{D_{i}}{2\pi \ii},\tau )\theta (-d\alpha _{i}-\varepsilon
\beta _{i}(c\tau +d),\tau )}\ee ^{-c D_{i} \varepsilon \beta _{i} },
\end{align*}
where in the last step we used the elliptic property
$$
\theta(z+\tau)=-\ee^{-\pi\ii(2z+\tau)} \theta(z,\tau).
$$

We can assume that the values $\beta _{i}$ are chosen so that $\beta
_{i}=m\cdot d_{i}$ for some generic $m\in M_\CC$.  This guarantees
$\sum \beta _{i} D_{i}$ is trivial in $H^{*} (X)$.  It is easy to see
that the arguments of Theorem~\ref{limit} still work, provided $\beta
_{i}\not =0$.  

Now we have 
\begin{align*}
f_{N,\deg} \left(\frac{a\tau +b}{c\tau +d} \right) (c\tau +d)^{-r} &= 
\lim _{\varepsilon \rightarrow 0}\int _{\hat{X}} \prod
_{i}\frac{(\frac{D_{i}}{2\pi \ii})\theta (\frac{D_{i}}{2\pi \ii}- d\alpha
_{i}-\varepsilon \beta _{i} (c\tau +d), \tau) \theta '(0,\tau
)}{\theta (\frac{D_{i}}{2\pi \ii},\tau )\theta (-d\alpha _{i}-\varepsilon
\beta _{i}(c\tau +d),\tau )}\\
&= f_{N,d\cdot \deg} (\tau ).
\end{align*}
\end{proof}

\subsection{}
To ease notation, we often suppress $\tau $ from $\theta $ and write
simply $\theta (z)$.  We introduce functions
$$
s_{a/l}^{(k)}(q) = (2\pi\ii)^{-k}
\left({\partial^k\over\partial z^k}\right)_{z=0}{\rm log}
\left(
{
z\theta(z+a/l)\theta'(0)
\over
\theta(z)\theta(a/l)
}
\right)
$$
for all $a=1,\dots ,l-1$ and all integers $k\geq 1$. Formally we will
also consider them for all integers $a$ not divisible by $l$ because
of the obvious periodicity.

Also we introduce
$$
r^{(k)}(q)=\left({\partial^k\over\partial z^k}\right)_{z=0}
{\rm log}\left({\theta(z)\over z\theta'(0)}\right)
$$
for all $k\geq 1$. We have $r^{(k)}=0$ for odd $k$, and
$s_{-a/l}^{(k)}=(-1)^ks_{a/l}^{(k)}.$

\begin{lemma}
\begin{enumerate}
\item For all $\left(\begin{array}{cc}a&b\\c&d\end{array}\right)\in
\Gamma_1(l)$,
$$
s_{a/l}^{(k)}\Bigl({a\tau+b\over c\tau+d}\Bigr)=(c\tau+d)^k s_{a/l}^{(k)}(\tau).
$$
\item For all $\left(\begin{array}{cc}a&b\\c&d\end{array}\right)\in
SL(2,\ZZ)$,
$$
r^{(k)}\Bigl({a\tau+b\over c\tau+d}\Bigr)=(c\tau+d)^k r^{(k)}(\tau),~k\geq 4
$$
$$
r^{(2)}\Bigl({a\tau+b\over c\tau+d}\Bigr)=(c\tau+d)^2 r^{(2)}(\tau)+
{c(c\tau+d)\over 2\pi\ii}.
$$
\end{enumerate}
\end{lemma}

\begin{proof}
This can be shown by arguments similar to those used to prove
Proposition \ref{modularity.of.tf}, and so we omit the details.
\end{proof}

\begin{proposition}\label{4.6}
Every element of ${\TTT}(l)$ is a polynomial with constant
coefficients of $s_{a/l}^{(k)}$ and $r^{(k)}$.  Moreover, the grading
of ${\TTT}(l)$ by the rank of the lattice coincides with the one
defined by the superscripts in brackets.
\end{proposition}

\begin{proof}
In the smooth case this result follows immediately from
Theorem~\ref{klem} and Taylor's formula.  Let $\alpha _{i} = \deg
(d_{i})$.  Then for each $i$ we have
$$
{(D_i/2\pi\ii)\theta(D_i/2\pi\ii-\alpha _{i})\theta'(0)
\over
\theta(D_i/2\pi\ii)\theta(-\alpha _{i})
}
=
\exp\Bigl (\sum_{k\geq 1}{(-1)^k\over k!}D_i^k s_{\alpha _{i}}^{(k)}\Bigr).
$$
Hence $f_{N,\deg}$ is polynomial in $s_{a/l}^{(k)}.$

Things are more complicated in the general case.
Let $\beta_i=\deg_1(d_i)$, where $\deg_{1}$ is defined in
Theorem~\ref{limit}. 
If $\alpha_i\notin\ZZ$, then
$${\log}
{(D_i/2\pi\ii)\theta(D_i/2\pi\ii-\alpha_i-\varepsilon \beta_i)\theta'(0)
\over
\theta(D_i/2\pi\ii)\theta(-\alpha_i-\varepsilon \beta_i)
}
=
{\rm log}
{(D_i/2\pi\ii-\varepsilon \beta_i)\theta(D_i/2\pi\ii-\alpha_i-\varepsilon
\beta_i)\theta'(0)
\over
\theta(D_i/2\pi\ii-\varepsilon \beta_i)\theta(-\alpha_i)
}
$$
$$
-
{\rm log}
{(-\varepsilon \beta_i)\theta(-\alpha_i-\varepsilon
\beta_i)\theta'(0)
\over
\theta(-\varepsilon \beta_i)\theta(-\alpha_i)
}
+
{\rm log}
{\theta(D_i/2\pi\ii-\varepsilon\beta_i)
\over
(D_i/2\pi\ii-\varepsilon\beta_i)\theta'(0)
}
-
{\rm log}
{\theta(D_i/2\pi\ii)
\over
(D_i/2\pi\ii)\theta'(0)
}
$$
$$
-
{\rm log}
{\theta(-\varepsilon\beta_i)
\over
(-\varepsilon\beta_i)\theta'(0)
}.
$$
This allows us to write 
$${(D_i/2\pi\ii)\theta(D_i/2\pi\ii-\alpha_i-\varepsilon \beta_i)\theta'(0)
\over
\theta(D_i/2\pi\ii)\theta(-\alpha_i-\varepsilon \beta_i)
}$$
as a power series in $D_i$ and $\varepsilon$, whose coefficients
are polynomials in $s_{\alpha_i}^{(k)}$ and $r^{(k)}$.

Now suppose 
$\alpha_i\in \ZZ $.  Then
$$
{\rm log}
{(D_i/2\pi\ii)\theta(D_i/2\pi\ii-\alpha_i-\varepsilon \beta_i)\theta'(0)
\over
\theta(D_i/2\pi\ii)\theta(-\alpha_i-\varepsilon \beta_i)
}
=
{\rm log}
{(D_i/2\pi\ii)\theta(D_i/2\pi\ii-\varepsilon \beta_i)\theta'(0)
\over
\theta(D_i/2\pi\ii)\theta(-\varepsilon \beta_i)
}
$$
$$
=
{\rm log}(D_i/2\pi\ii-\varepsilon \beta_i)
-{\rm log}(-\varepsilon \beta_i)
+{\rm log}
{\theta(D_i/2\pi\ii-\varepsilon\beta_i)
\over
(D_i/2\pi\ii-\varepsilon\beta_i)\theta'(0)
} 
-
{\rm log}
{\theta(D_i/2\pi\ii)
\over
(D_i/2\pi\ii)\theta'(0)
}
$$
$$
-
{\rm log}
{\theta(-\varepsilon\beta_i)
\over
(-\varepsilon\beta_i)\theta'(0)
}.
$$
Thus we can write
$$
\prod_i
{(D_i/2\pi\ii)\theta(D_i/2\pi\ii-\alpha _{i}-\varepsilon \beta _{i})\theta'(0)
\over
\theta(D_i/2\pi\ii)\theta(-\alpha _{i}-\varepsilon \beta _{i})
}
$$
as $\varepsilon^{-l}$ times a power series of $D$ and $\varepsilon$
whose coefficients are polynomial of $s_{a/l}^{(k)}$ and $r^{(k)}$. 
Since the limit as $\varepsilon\to 0$ exists, after we integrate we
find no negative degrees of
$\varepsilon$. 

Finally, the statement about grading is proved by looking at the total
degree in $D$ and $\varepsilon$ and observing that at the end we put
$\varepsilon=0$.
\end{proof}

\subsection{}
We will now show that the $s_{a/l}^{(k)}$ can be expressed
as polynomials in $s_{a/l}^{(1)}$ with minor exceptions.

\begin{lemma}
Let $N$ be any natural number and let $a_1,\dots ,a_{N+1}$ be
a set of nonzero residues $\mod l$ so that
$a_1+\cdots+a_{N+1}= 0 \mod l $. Then
the coefficient of $t^N$ in 
$$
\exp\Bigl(\sum_{k\geq 0}(t^k/k!)\bigl(\sum_{j=1}^{N+1} s_{a_j/l}^{(k)}\bigr)\Bigr)
$$   
is zero.
\label{relation}
\end{lemma}
  
\begin{proof}
From the definition of $s_{a/l}^{(n)}$ and Taylor's formula,
\[
\exp\Bigl(\sum_{k\geq 0}(t^k/k!)\bigl(\sum_{j=1}^{N+1}
s_{a_j/l}^{(k)}\bigr)\Bigr)
=\frac{t^{N+1}}{(2\pi \ii)^{N+1} }\prod_{j=1}^{N+1}
{\theta(t/2\pi \ii +a_j/l)\theta'(0) \over \theta(t/2\pi \ii )\theta(a_j/l)}.
\]

Let 
\[
F(t,\tau)=\prod_{j=1}^{N+1} {\theta(t+a_j/l)\theta'(0)  
\over
\theta(t)\theta(a_j/l)}.
\]
We want to show the residue of $F$ at $t=0$ is zero. But notice that
$F$ is a doubly-periodic function due to the condition on $\sum
a_j$. Therefore, the sum of its residues in a fundamental domain is
zero.  But $F$ has a unique pole at $t=0$, which finishes the
proof.
\end{proof}

\begin{proposition}
\begin{enumerate}
\item If $l\geq 5$ then all $s_{a/l}^{(k)}$
are polynomials of $s_{a/l}^{(1)}$.
\item If $l=4$ then all $s_{a/4}^{(k)}$ are polynomials of 
$s_{1/4}^{(1)}$ and $s_{1/4}^{(2)}$.
\item If $l=3$ then all $s_{a/3}^{(k)}$ are polynomials of 
$s_{1/3}^{(1)}$ and $s_{1/3}^{(3)}$.
\item If $l=2$ then all $s_{1/2}^{(k)}$ are polynomials of 
$s_{1/2}^{(2)}$ and $s_{1/2}^{(4)}$.
\end{enumerate}
\label{saone}
\end{proposition}

\begin{proof}
{\em Case $l\geq 5$.}
We need to show that for a fixed $N\geq 2$ all $s_{a/l}^{(N)}$
can be expressed in terms of $s_{a/l}^{(\leq N-1)}$.
Denote the ring generated by $s_{a/l}^{(\leq N-1)}$ by
$R_{N-1}$. Because of Lemma \ref{relation}, for every 
set of nonzero residues $a_1,\dots ,a_{N+1}$ with $\sum_j a_j = 0
\mod l$, we have
$$\sum_j s_{a/l}^{(N)} = 0 \mod R_{N-1}.$$

If $N\geq 3$, then this implies that for any four residues
$a,b,c,d$ with $a+b=c+d$, we have 
$$s_{a/l}^{(N)} + s_{b/l}^{(N)} = s_{c/l}^{(N)} + s_{d/l}^{(N)}
\mod R_{N-1}.
$$
This implies 
$$s_{a/l}^{(N)} =  (2-a) s_{1/l}^{(N)} + (a-1)s_{2/l}^{(N)} \mod
R_{N-1},\quad a=1,\dots ,l-1.$$
We now go back to the original relation and see that if 
$\sum_{j=1}^{N+1} a_j = kl$, then 
$$(2(N+1)-kl)s_{1/l}^{(N)}+(kl-(N+1))s_{2/l}^{(N)}=0 \mod R_{N-1}.$$
Since $l\geq 5$ we can find two sets of $a_j$ with 
different $k$, which shows that $s_{1/l}^{(N)},s_{2/l}^{(N)}\in
R_{N-1}$.

If $N=2$, we have 
$$s_{a/l}^{(2)}+ s_{b/l}^{(2)}+ s_{c/l}^{(2)}=0\mod R_1,
$$
if $a+b+c = 0\mod l$. Together with the symmetry 
$s_{a/l}^{(2)}=s_{(l-a)/l}^{(2)}$, this implies 
$$s_{(a+1)/l}^{(2)}=-s_{a/l}^{(2)}-s_{1/l}^{(2)}\mod R_1$$
if $a,a+1\neq 0\mod l$. Therefore,
$$s_{a/l}^{(2)} = ({-1-3(-1)^{a}\over2})s_{1/l}^{(2)}\mod R_1$$
for $a=1,\dots ,l-2$.
It remains to consider $a=b=2, c=l-4$ to show that all $s_{a/l}^{(2)}$
lie in $R_1$.

{\em Case $l=4$.} The only difference is that for $N=2$ one can no 
longer use $a=b=2,c=l-4$.

{\em Case $l=3$.} At $N=3$ we can no longer find two sets of four
$a_j$ with different $k$. On the other hand at $N=2$ we easily
conclude that $s_{1/3}^{(2)}\in R_1$, because of $1+1+1=0\mod l$.

{\em Case $l=2$.}  By symmetry, all we have is $s_{1/2}^{(2k)}$ for
$k\geq 1$. We notice that $s_{1/2}^{(2k)}$ is proportional to
$G^*_{2k}$ of \cite{Zagier.elliptic}, and is therefore a modular
form. It is also known (see \cite{Landweber.elliptic}) that the ring
of modular forms for $\Gamma_0(2)$ is freely generated by an element
of degree $2$ and an element of degree $4$. It remains to observe that
$(s_{1/2}^{(2)})^2$ is not proportional to $s_{1/2}^{(4)}$.
\end{proof}

\begin{lemma}
For all $l$ and $n\geq 4$, any $r^{(n)}$ is a polynomial in
$s_{a/l}^{(k)}$. 
\label{efroms}
\end{lemma}

\begin{proof}
Consider the function
$$f(z,\tau)={\partial^2\over\partial z^2}\log(\theta(z,\tau)).$$
This is easily seen to be a doubly-periodic function. Its only
pole in the fundamental domain is at the origin, and the Laurent expansion
around it is
$$f(z,\tau) = {(-1)\over z^2} + 
\sum_{k\geq 0} r^{(k+2)}(\tau){z^k\over k!}.$$
On the other hand, consider the function 
$$g(z,\tau) = 
\left(
{\theta(z+1/l,\tau)\theta'(0,\tau)
\over
\theta(1/l,\tau)\theta(z,\tau)
}
\right)
\cdot
\left(
{\theta(z-1/l,\tau)\theta'(0,\tau)
\over
\theta(-1/l,\tau)\theta(z,\tau)
}
\right)
$$
The function $g$ is even, so it has no residue at $z=0$.
One easily sees that its Laurent expansion is 
$$g(z)={1\over z^2} + \varphi (z),$$
where $\varphi (z)$ is holomorphic.
Thus, $f+g$ is an even doubly-periodic holomorphic
function, which implies that it is constant in $z$. 
As a result, the coefficients of its Laurent expansion
at $z^n,~n\geq 2$ are zero. It remains to notice that all
coefficients of $g$ at $z=0$ are polynomials in $s_{1/l}^{(k)}$.
\end{proof}

\begin{theorem}
The ring ${\TTT}(l)$ is a subring of the ring of modular
forms for $\Gamma_1(l)$. 
\begin{enumerate}
\item If $l\geq 5$, then it is generated by 
$s_{a/l}^{(1)}$. 
\item If $l=4$, then it is generated
by $s_{1/4}^{(1)}$ and  $s_{1/4}^{(2)}$. 
\item If $l=3$, then it is generated by $s_{1/3}^{(1)}$ and
$s_{1/3}^{(3)}$.
\item If $l=2$, then it is generated by  $s_{1/2}^{(2)}$ and $s_{1/2}^{(4)}$.
\end{enumerate}
\label{main}
\end{theorem}

\begin{proof}
Because of Lemma \ref{efroms}, we can express every element of
${\TTT}(l)$ as a polynomial in $r^{(2)}$ and $s_{a/l}^{(k)}$.  The
automorphic properties allow us to conclude that $r^{(2)}$ is
absent from these polynomial expressions. Together with
Proposition~\ref{saone}, this implies that ${\TTT}(l)$ is contained inside
the ring generated by $s_{a/l}^{(1)}$ if $l\geq 5$, and by
appropriate $s_{a/l}^{(k)}$ for $l=2,3,4$.

To prove the opposite inclusion, we first notice that ${\TTT}(l)$
contains all $s_{a/l}^{(1)}$. Indeed, consider the degree function on
the one-dimensional lattice $N=\ZZ$ that equals $a/l$ on $n=1$
and $n=-1$.  There is only one possible complete fan, and the
corresponding toric variety is $\Proj^1$.  By Theorem~\ref{klem}, we get
$$
f_{N,\deg}= -2s^{(1)}_{a/l}.
$$ 
Similarly, by using $\Proj^k$ we conclude
that $s_{a/l}^{(k)}\in {\TTT}(l)$, which helps us at levels
less than $5$.

To finish the proof, we must show that toric forms have the expected
behavior at all cusps.  Let $\alpha \in \QQ $.  Here is the argument
for $s_{\alpha }^{(1)}$. It is easy to see that
\begin{equation}\label{what.sa.is}
s_{\alpha }^{(1)}(\tau)={\theta'(\alpha ,\tau)\over 2\pi \ii
\theta(\alpha ,\tau)}.
\end{equation}
We need to show that for any
$\left(\begin{array}{cc}a&b\\c&d\end{array}\right)\in \Gamma(1)$, the
function
$$s_{\alpha }^{(1)}\Bigl({a\tau+b\over c\tau+d}\Bigr) (c\tau+d)^{-1}$$ 
is bounded in any neighborhood of $\ii\infty$. 
Using the transformation properties of $\theta$, we can deduce that
$$
s_{\alpha }^{(1)}\Bigl({a\tau+b\over c\tau+d}\Bigr) (c\tau+d)^{-1}=
{\theta'(c\alpha \tau+d\alpha ,\tau)\over 2\pi\ii \theta(c\alpha \tau+d\alpha ,\tau)}+{cd\alpha \over
c\tau+d}.
$$
Using elliptic properties of $\theta$, we can reduce $c\alpha
\tau+d\alpha $ to $\lambda\tau+\nu$ where $0\leq\lambda,\nu< 1$ and
$(\lambda,\nu)\neq(0,0)$. This introduces an extra summand that is
irrelevant, since all we want to show is boundedness. We will closely
examine the $q$-expansion of $\theta(\lambda \tau+\nu,\tau)$.

From the definition of $\theta$ \eqref{theta.def}, we have
$$\theta(\lambda\tau+\nu,\tau)={1\over i}\sum_{n\in\ZZ}(-1)^n
\ee^{\pi\ii\tau((n+{1\over 2})^2+\lambda(2n+1))} \ee^{\pi \ii \nu
(2n+1)}.$$ 
Since $0\leq \lambda <1$, the smallest (rational) power
of $q$ can occur at $n=-1$ if $\lambda\neq 0$ or at
$n=-1,0$ if $\lambda=0$. 

In the first case, it is clear that the coefficient is nonzero.  Since
differentiation by $z$ contributes a nonzero factor of $\pi\ii(2n+1)$,
the ratio in the formula above starts with $q^0$, and has nonnegative
powers of $q$.

In the second case, two terms that contribute to the smallest power of
$q$ have nonzero coefficients whose ratio is $(-1)\ee^{2\pi \ii
\nu}$.  Hence there is no cancellation.  As a result, the ratio
again has only nonnegative powers of $q$.

A similar argument works for all $s_{a/l}^{(k)}$, and so this
completes the proof.
\end{proof}

\begin{remark}
It is easy to see that all elliptic genera of Hirzebruch that are
modular with respect to $\Gamma_1(l)$ (see \cite{hirz}) take values in
${\TTT}(l)$ for appropriate $l$. In general, however, they form a
subspace, because in their definition one only uses $s_{\alpha}^{(k)}$ for
a fixed $\alpha$. Nevertheless, ${\TTT}(l)$ could be alternatively
defined as the space spanned by all products of all Hirzebruch elliptic
genera that are modular with respect to $\Gamma_1(l)$.
\end{remark}

\begin{remark}
Even though we will not need it for our arguments, it is worth
mentioning that $r^{(k)}=-2E_k$ for all (even) $k$ where $E_k$ is the
standard Eisenstein series.  To prove it, use the product formula for
$\theta(z,\tau)$, take its logarithm, differentiate once with respect
to $z$, expand the result as a geometric series and differentiate some
more. So the ring generated by $r^{(k)}$ is the ring of almost-modular
forms, see \cite {Zagier.elliptic}.  \label{almostmod}
\end{remark}

\begin{remark}\label{conjecture}
It is natural to ask whether ${\TTT}(l) = \MMM (\Gamma _{1} (l))$.
For example, it is not hard to show that $\TTT (l) = \MMM (\Gamma
_{1} (l))$ for $l=2$, $3$, and $4$.  In general, however, the answer is no.
In particular, we can never construct any cusp forms of weight one this
way.

We have investigated this question in more detail for the case of
weight $2$ and prime level.  Specifically, we have considered the
$\Gamma_0(l)/\Gamma _{1} (l)$-invariant part of ${\TTT}_{2}(l)$, which
consists of modular forms for $\Gamma _{0} (l)$.  Extensive computer
calculations suggest that $({\TTT}_{2}(l))^{\Gamma_0(l)/\Gamma _{1}
(l)}$ coincides with the span of the Eisenstein series and the
eigenforms of analytic rank zero.  We will address this conjecture in
more detail in \cite{nextpaper}.
\end{remark}

%
%

\section{Hecke stability}\label{hecke.sect}

\subsection{}
To conclude the paper, we will show that toric forms of level $l$ are
well-behaved with respect to the Hecke operators, the involution
operator $w_{l}$, and Atkin-Lehner lifting.
For background, we refer to Lang \cite{Lang}.

We begin with the Hecke operators.  These are endomorphisms
of $\MMM_{r} (\Gamma _{1} (l))\otimes \CC $ generated as follows.  For 
a prime $p$ with $(p,l)=1$, the operator $T_{p}$ acts by
\begin{equation*}
(f\bigm|T_p)(\tau) = p^{-1}\sum_{i=0}^{p-1}f((\tau+i)/p) +
p^{r-1} (f\bigm|\varepsilon_p)(p\tau),
\end{equation*}
Here $f\bigm|\varepsilon_p$ denotes the action of an element of
$\Gamma_0(l)/\Gamma_1(l) \cong (\ZZ /l\ZZ )^{\times }$ given by a
matrix
$$
\left(\begin{array}{cc}a&b\\c&d\end{array}\right)=
\left(\begin{array}{cc}p^{-1}&b\\0&p\end{array}\right)
\mod l.
$$

For a prime $p$ with $(p,l)>1$, the operator $U_{p}$ acts by
\begin{equation*}
(f\bigm|U_p)(\tau) = p^{-1}\sum_{i=0}^{p-1}f((\tau+i)/p)
\end{equation*}
In terms of $q$-expansions, if $f = \sum _{n}a_{n}q^{n}$, we have
\begin{equation*}
f\bigm|U_{p} = \sum _{p|n} a_{n}q^{n/p}.
\end{equation*}
If we define an operator $V_{p}$ by
\begin{equation}\label{vop}
f\bigm|V_{p} = \sum _{n}a_{n}q^{np},
\end{equation}
then we have for $(p,l)=1$
\begin{equation}\label{qexp1}
T_{p} = U_{p} + p^{r-1}\varepsilon _{p}V_{p}.
\end{equation}

It is easy to see how $\varepsilon _{p}$ acts on the toric
forms.
\begin{lemma}
Let $f_{N,\deg}\in \TTT (l)$.  If $(p,l)=1$, then
$$
f_{N,\deg}\bigm|\varepsilon_p = f_{N,p\deg}.
$$ 
\label{diam}
\end{lemma}

\begin{proof}
This follows from Proposition~\ref{modularity.of.tf}.
\end{proof}

\begin{theorem}
The space $\TTT_{r}(l)$ is stable under the action of the Hecke 
operators.
\label{hecke}
\end{theorem}

\begin{proof}
First we assume $(p,l)=1$ and consider the operator $T_{p}$.  Let
$f_{N,\deg}\in \TTT _{r} (l)$, and assume $\deg$ is linear with respect to
$\Sigma $.

Using Lemma~\ref{diam} and \eqref{qexp1}, we find
\begin{align*}
f_{N,\deg}\bigm|T_p&=\sum_{m\in M} \sum_{C \in \Sigma}(-1)^{\codim C }
\ac{\sum_{n\in C \cap N} \delta_{m\cdot n}^{0\mod p\ZZ } 
q^{m\cdot n/p}\ee^{2\pi\ii \deg(n)}}\\
&\qquad + p^{r-1}\sum_{m\in M}\sum_{C \in
\Sigma}(-1)^{\codim C }\ac{\sum_{n\in C \cap N} q^{m\cdot np}\ee^{2\pi\ii p \deg(n)}}\\
&=\sum_{m\in M} \sum_{C \in \Sigma}(-1)^{\codim C } \ac{\sum_{n\in C
\cap N} (\delta_{m\cdot n}^{0\mod p\ZZ}
q^{m\cdot n/p}\ee^{2\pi\ii \deg(n)}\\
&\qquad + p^{r-1}\delta_{m}^{0\mod pM} q^{m\cdot n}\ee^{2\pi\ii p \deg(n)})}.
\end{align*}
Here $\delta ^{\alpha }_{\beta }$ takes the value $1$ if $\beta $ is
the residue $\alpha $, and is $0$ otherwise.

We will compare this expression with 
\[
\sum _{S} f_{S,p\deg},
\]
where the sum is taken over lattices $S$ satisfying $N\subset S\subset
\frac{1}{p}N$ and $[S:N] = p^{r-1}$.  The dual of each $S$ lies in
$M$, and we have
\begin{align}
\sum_{S}f_{S,p\deg}
&=
\sum_{m\in M} \sum_{C \in \Sigma}(-1)^{\codim C }
\ac{\sum_{S}\sum_{n\in C \cap S} \delta_{m\cdot S}^{0\mod\ZZ}
q^{m\cdot n}\ee^{2\pi\ii p \deg (n)}}\\
&=\sum_{m\in M} \sum_{C \in \Sigma}(-1)^{\codim C }
\ac{\sum_{R}\sum_{n\in C \cap R} \delta_{m\cdot R}^{0\mod p\ZZ}
q^{m\cdot n/p}\ee^{2\pi\ii \deg (n)}},\label{rline}
\end{align}
Here the sum over $R$ in \eqref{rline} is taken over lattices with
$N\supset R \supset pN$ and $[N:R]=p$.

Now consider the contributions to the $R$-sum by different $m\in M$.
If $m\not \in pM$, then there is only one $R$ for which $m\cdot R$ is
always divisible by $p$: the set of all $n$ with $m\cdot n = 0
\mod p$. So for every $m\notin pM$, each cone $C \in \Sigma$
contributes
$$
\ac{\sum_{n\in C \cap N}\delta_{m\cdot n}^{0\mod p\ZZ}q^{m\cdot
n/p}\ee^{2\pi\ii \deg (n)}}.
$$

On the other hand, if $m\in pM$, then the sum \eqref{rline} is taken
over all $R$.  Letting $\mu (n) = \sharp\{ R\mid n\in
R\}$, we have for fixed $m$ and $C $
\begin{align*}
&\ac{\sum_{R}\sum_{n\in C \cap R}
q^{m\cdot n/p}\ee^{2\pi\ii \deg (n)}}=
\ac{\sum_{n\in C \cap N}\mu (n) q^{m\cdot n/p}\ee^{2\pi\ii \deg (n)}}\\
&={p^{r-1}-1\over p-1}\ac{\sum_{n\in C \cap N}q^{m\cdot n/p}\ee^{2\pi\ii
\deg(n)}}
+p^{r-1}\ac{\sum_{n\in C \cap pN}q^{m\cdot n/p}\ee^{2\pi\ii \deg (n)}}\\
&={p^{r-1}-1\over p-1}\ac{\sum_{n\in C \cap N}q^{m\cdot n/p}\ee^{2\pi\ii
\deg(n)}} +p^{r-1}\ac{\sum_{n\in C \cap N}q^{m\cdot n}\ee^{2\pi\ii
p\deg (n)}}.
\end{align*}

Collecting these facts, we see
\begin{align*}
f_{N,\deg}\bigm|T_p-\sum_{S}f_{S,p\deg}(q)
&=-{p^{r-1}-p\over p-1}\sum_{m\in pM}\sum_{C \in \Sigma}(-1)^{{\rm
codim}C }\ac{\sum_{n\in C }q^{m\cdot n/p}\ee^{2\pi\ii\deg(n)}}\\
&=-{p^{r-1}-p\over p-1}\sum_{m\in M}\sum_{C \in \Sigma}(-1)^{{\rm
codim}C }\ac{\sum_{n\in C }q^{m\cdot n}\ee^{2\pi\ii\deg(n)}}\\
&={p-p^{r-1}\over p-1}f_{N,\deg}.
\end{align*}
This shows that $f_{N,\deg}\bigm| T_{p}$ is a linear combination of
toric forms.

Now let $(p,l)>1$, and consider the operator $U_{p}$.  To show that
$\TTT _{r} (l)$ is stable under $U_{p}$, we use a technique similar to 
the proof of Theorem~\ref{limit}.  Let $\varepsilon $ be a parameter,
and let $\deg'$ be a generic degree function with respect to
$\Sigma $.  By the argument above, we have 
\[
f_{N,\deg + \varepsilon \deg'}\bigm| T_{p} = \sum _{S_{i}} \alpha
_{i}f_{S_{i},\deg_{i} + \varepsilon \deg_{i}'},
\]
where the sum is taken over some collection of lattices $\{S_{i} \}$
with $\deg_{i} (S_{i})\subset \frac{1}{l}\ZZ $ and $\deg_{i}'$ generic
for all $i$.

By \eqref{vop}, we have
\begin{equation}\label{differ}
f_{N,\deg+\varepsilon \deg'}\bigm|U_{p} = f_{N,\deg+\varepsilon \deg'}
\bigm| T_{p} - p^{r-1}f_{N,p\deg + \varepsilon p \deg '}\bigm| V_{p}.
\end{equation}
Let $m_{1},\dots ,m_{k}\subset M$ be a complete set of residues for
$M/pM$.  For any $m\in M$, we write $[m/p]$ for the function $n
\mapsto (m\cdot n)/p$.

\begin{lemma}\label{vlemma}
We have
\[
\sum _{k} f_{N,\deg + \varepsilon \deg' + [m_{k}/p]} = p^{r}
f_{N,p\deg+\varepsilon p\deg'}\bigm|V_{p}.
\]
\end{lemma}

\begin{proof}[Proof of Lemma \ref{vlemma}]
We have 
\begin{multline*}
\sum _{k} f_{N,\deg + \varepsilon \deg' + [m_{k}/p]} \\
= \sum _{m\in
M}\sum _{C\in \Sigma } (-1)^{\codim C} \ac{\sum _{n\in C} q^{m\cdot
n}\sum _{k} \ee ^{2\pi \ii (\deg + \varepsilon \deg ' + [m_{k}/p]) (n)}}. 
\end{multline*}
Using 
\[
\sum _{k} \ee ^{2\pi \ii (m_{k}\cdot n)/p} = \begin{cases}
0&n\not \in pN\\
p^{r}&n\in pN
\end{cases}
\]
this becomes 
\[
p^{r} \sum _{m\in M}\sum _{C\in \Sigma } (-1)^{\codim C} \ac{\sum
_{n\in C\cap pN} q^{m\cdot n} \ee ^{2\pi \ii (\deg + \varepsilon \deg
') (n)}} = p^{r} f_{N,p\deg+\varepsilon p\deg'}\bigm|V_{p}.
\]
\end{proof}

Lemma~\ref{vlemma} and \eqref{differ} imply 
\begin{equation}\label{uop}
f_{N,\deg + \varepsilon \deg'}\bigm|U_{p} = \sum_{i} \alpha _{i}
f_{S_{i},\deg_{i}+\varepsilon \deg_{i}'},
\end{equation}
for some collection of lattices $S_{i}$ and degree functions
$\deg_{i}, \deg'_{i}$.  We can desingularize $\Sigma $ separately with respect to each $S_{i}$.  Then by Theorem~\ref{klem}, 
\eqref{uop} equals
\begin{equation*}
\sum _{i} \int _{X_{i}} \Phi _{i} (\varepsilon ),
\end{equation*}
where each $X_{i}$ is a nonsingular toric variety, and $\Phi _{i}
(\varepsilon )$ is an $\varepsilon $-family of rational functions in
$\theta $ and $\theta '$.  

Since the action of $U_{p}$ commutes with the $\varepsilon $-limit, we
see that
\begin{align*}
f_{N,\deg }\bigm|  U_{p} &= \lim_{\varepsilon \rightarrow 0} \sum _{i}
\int _{X_{i}} \Phi _{i} (\varepsilon )\\
&=\lim_{\varepsilon \rightarrow 0} P (\varepsilon ).
\end{align*}
Here $P (\varepsilon )$ is a Laurent series in $\varepsilon $ with
coefficients that are polynomials in $s_{a/l}^{(1)}$ and $r^{(2)}$
($l\geq 5$), and by appropriate $s^{(k)}$ for other $l$, as in
Proposition~\ref{4.6}.  This limit exists and converges to a modular
form.  Hence $r^{(2) }$ doesn't appear in the final expression, and by
Theorem~\ref{main}
\[
f_{N,\deg }\bigm|  U_{p} \in \TTT _{r} (l).
\]
\end{proof}

\subsection{}
As a by-product of the above argument, we get that toric forms are
stable under lifting of oldforms.

\begin{corollary}\label{atkin}
Let $f (q) = f (\tau )$ be a toric form of level $l$.  Then $f (p\tau
)$ is a sum of toric forms of level $pl$.
\end{corollary}

\begin{proof}
As in 
Lemma~\ref{vlemma},
\[
f_{N,\deg} (p\tau ) = p^{-r}\sum _{k} f_{N,\deg/p + [m_{k}/p]} (\tau ).
\]
\end{proof}

\begin{remark}
At this point we do not know any direct proof of Hecke stability
in terms of $s_{a/l}^{(k)}$. Because Hecke operators are generally
not compatible with multiplication, we find it fascinating
that there exists a Hecke stable finitely generated subring of the ring
of modular forms.
\end{remark}

\subsection{}
Finally, we consider the Fricke involution.  Let $W_{l}$ be the matrix
\[
\left( 
\begin{array}{cc}
0&-1\\
l&0
\end{array}\right)
\]

\begin{theorem}\label{fricke}
The space $\TTT (l)$ is stable under the action of $W_{l}$.
\end{theorem}

\begin{proof}
By Remark~\ref{conjecture}, we need only consider $l\geq 5$, and for these
levels, it suffices to verify the statement for $s_{a/l}^{(1)}$.

We have
\[
\log \theta (z/l\tau + a,-1/l\tau ) = f (\tau ) + \frac{\pi \ii}
{l\tau} (z+a\tau )^{2} + \log\theta (z+a\tau ,l\tau ). 
\]
After differentiating with respect to $z$ and evaluating at $z=0$, we
obtain
\begin{equation}\label{sa.exp}
\frac{1}{l\tau }s_{a/l}^{1}\left(-\frac{1}{l\tau } \right) =
\frac{a}{l} + \frac{\theta' (a\tau ,l\tau ) }{2\pi \ii \theta (a\tau
,l\tau )}. 
\end{equation}

Now differentiate \eqref{jf} with respect to $z$ to obtain
\[
\frac{1}{2\pi \ii }\frac{\partial }{\partial z}\log\theta = \frac{1}{2} - \sum _{n\geq 1}\frac{\ee
^{2\pi\ii z}q^{n}}{1-\ee ^{2\pi \ii z}q^{n}} + \sum
_{n\geq 0}\frac{\ee ^{-2\pi\ii z}q^{n-1}}{1-\ee ^{-2\pi \ii
z}q^{n-1}}.
\]
The right-hand side of \eqref{sa.exp} becomes
\begin{equation}\label{eq2}
C- \sum _{n\geq 1}\frac{q^{a+nl}}{1-q^{a+nl}} + \sum
_{n\geq 0}\frac{q^{-a+nl}}{1-q^{-a+nl}},
\end{equation}
where $C$ is a constant.

In \eqref{eq2} we can take the first sum over $n\geq 0$ and the second
over $n\geq 1$, and absorb the correction in the constant.  After
simplifying, we find that \eqref{eq2} becomes
\begin{equation*}
C - \sum _{d\geq 1}q^{d} \bigl( \sum _{a+nl|d}-1 + \sum
_{-a+nl|d}+1\bigr) = C - \sum _{d\geq 1}q^{d}\sum _{k|d} (\delta _{k}^{a\mod l} - \delta
_{k}^{-a\mod l}), 
\end{equation*}
where $\delta $ is defined as in the proof of Theorem~\ref{hecke}.

Now it is not hard to show, using \eqref{what.sa.is}, that
\[
s_{b/l}^{(1)} = C_{1}(b) - \sum _{d}q^{d}\sum _{k|d} (\ee ^{2\pi \ii
kb/l}-\ee ^{-2\pi \ii kb/l}).
\]
Convolving with roots of unity, we obtain
\begin{align*}
\sum _{j=1}^{l-1}s_{j/l}^{(1)}\ee ^{-2\pi \ii ja/l} &= C_{2}- \sum _{d}q^{d}\sum _{k|d}\bigl (\sum _{j=1}^{N-1}\ee ^{-2\pi \ii ja/l} (\ee ^{2\pi \ii kj/l}-\ee ^{-2\pi \ii kj/l})\bigr)\\
&=C_{2}-l\sum _{d}q^{d}\sum _{k|d}\left (\delta _{k}^{a\mod l} - \delta
_{k}^{-a\mod l}\right).
\end{align*}
Hence after subtracting a linear combination of $s_{b/l}^{(1)}$ from
$s_{a/l}^{(1)}\bigm|W_{l}$, we obtain a constant.  Since this
difference is modular, the constant must vanish, and the proof is complete.
\end{proof}

\begin{remark}\label{lang.rmk}
Notice that \eqref{eq2} implies that $s_{a/l}^{(1)}\bigm|W_{l}$ is the
same as the Hecke-Eisenstein form $G_{1,a}$ appearing in \cite[Ch. 15,
\S1]{Lang}.
\end{remark}

%
%

\section{Appendix: background on toric geometry}\label{appendix}

\subsection{}
In this section we collect basic facts from toric geometry.  For more
details, the reader may consult \cite{fulton}.

\subsection{}
Let $N$ be a lattice, and let $N_{\RR} = N\otimes \RR$.  A subset
$C\subset N_{\RR}$ is a \emph{cone} if $C$ is closed under homotheties
and contains no line.  

Let $M=\Hom _{\ZZ } (N,\ZZ )$ be the dual lattice, and let $m\cdot n$
be the natural pairing.  A cone $C$ is \emph{rational polyhedral} if
there exist $m_{1},\dots ,m_{s}\in M$ such that 
\begin{equation}\label{def.of.rat}
C = \bigcap _{i} \bigl\{x\in N_{\RR}\bigm| m_{i}\cdot x \geq 0\bigr\}.
\end{equation}
A \emph{face} of $C$ is the subset of $C$ obtained by making some of
the inequalities in \eqref{def.of.rat} equalities.
A rational polyhedral cone $C$ is \emph{simplicial} if 
\[
C = \sum _{i=1}^{k} \RR _{\geq 0} n_{i},\quad n_{i}\in N,
\]
where $k$ is the dimension of the subspace generated by $C$.

\subsection{}
Let $\Sigma $ be a set of rational polyhedral cones.  Then $\Sigma $
is called a \emph{fan} if 
\begin{enumerate}
\item each face of a cone in $\Sigma $ is also in $\Sigma $, and 
\item the intersection of any two cones in $\Sigma $ is a face of each.
\end{enumerate}
A fan is \emph{complete} if 
\[
\bigcup_{C\in \Sigma } C = N_{\RR},
\]
and is simplicial if all $C\in \Sigma $ are simplicial.

\subsection{}
Let $\Sigma $ be a fan.  We can associate a \emph{toric variety}
$X_{\Sigma }$ to $\Sigma $ as follows.

If $C\subset
N_{\RR}$, let $C^{*}\subset M_{\RR}$ be the dual cone.  Then we have a
$\CC$-algebra $S_{C}$ defined by
\[
S_{C} = \CC [C^{*}\cap M],
\]
and we let $U_{C}$ be the affine variety $\Spec S_{C}$.  The variety
$U_{C}$ is the \emph{toric chart} associated to $C$. 

We can glue these affine varieties together using the combinatorics of
$\Sigma $.  Let $C_{1},C_{2}\in \Sigma $ satisfy $C_{1}\cap C_{2}\not
= \varnothing $.  Then we can identify $U_{C_{1}\cap C_{2}}$ with a
principal open subvariety of $U_{C_{1}}$ and of $U_{C_{2}}$.  These
identifications are compatible as we range over all the cones of
$\Sigma $, and so after gluing we obtain a variety $X_{\Sigma }$.  One
can show that $X_{\Sigma }$ is separated, and has an open set
isomorphic to an algebraic torus $T$.  Moreover, $X_{\Sigma }$ is
complete if and only if $\Sigma $ is complete, and $X_{\Sigma }$ is
nonsingular if and only if each cone of $\Sigma $ is generated by part
of a basis of $N$.

\subsection{}
Let $\{d_{i} \}\subset \Sigma $ be the set of one-dimensional cones.
Associated to each $d_{i}$ is a $T$-invariant divisor $D_{i}\subset
X_{\Sigma }$.  One can show that if $X_{\Sigma }$ is nonsingular, then
the cohomology classes of the $D_{i}$ generate $H^{*} (X_{\Sigma },\CC )$.
The linear relations among these classes are of the form 
\[
\sum _{i} ( m\cdot d_{i})D_{i} = 0,
\] 
where $m\in M_{\CC }$.
The total Chern class of $X$ is given by 
\[
c (TX) = \prod _{i} (1+D_{i}).
\]
From this it follows that the Todd class is given by 
\[
\Td (TX) = \prod_{i} \frac{D_{i}}{1-\ee ^{-D_{i}}}.
\]

%
%

\bibliographystyle{amsplain}
\bibliography{toric}

%
%

\end{document}